\documentclass[10pt,twoside]{article}

\author{\bf {D.D.\ Poro\c sniuc}\\
\it \small Department of Mathematics, National College "Mihai Eminescu",\\
\it \small RO-710096, Boto\c sani, Rom\^ ania.\\
\it \small e-mail: dporosniuc@yahoo.com\\
~~~~~~~~\it \small danielporosniuc@lme.ro}
\date{}
\title{\bf A class of locally symmetric \\K\"ahler Einstein structures on the cotangent
bundle of a space form}

\begin{document}

\maketitle \pagestyle{myheadings}\markboth {\small D.D. Poro\c
sniuc: A class of locally symmetric K\"ahler Einstein
structures...\hspace{1.2cm}}{\small \hspace{1.2cm}D.D. Poro\c
sniuc: A class of locally symmetric K\"ahler Einstein
structures...}

{\bf Abstract.} We obtain a class of locally symmetric K\"ahler
Einstein structures on the cotangent bundle of a Riemannian
manifold of negative sectional curvature. Similar results are
obtained in the case of a Riemannian manifold of positive
sectional curvature. The obtained class of K\"ahler Einstein
structures depends on one essential parameter and cannot have
constant holomorphic sectional curvature.

\vskip 0.5cm



\vskip 10mm {\large \bf{Introduction}} \vskip3mm

 In the study of the differential geometry of the cotangent bundle
$T^*M$ of a Riemannian manifold $(M,g)$ one uses several
Riemannian and semi-Riemannian metrics, induced from the
Riemannian metric $g$ on $M$. Next, one can get from $g$ some
natural almost complex structures on $T^*M$. The study of the
almost Hermitian structures induced from $g$ on $T^*M$ is an
interesting problem in the differential geometry of the cotangent
bundle.

 In \cite{OprPor2} the authors have obtained a class of natural
K\"ahler Einstein structures $(G,J)$ of diagonal type induced on
$T^*M$ from the Riemannian metric $g$. The obtained K\"ahler
structures on $T^*M$ depend on two essential parameters $a_1$ and
$\lambda$, which are smooth functions depending on the energy
density $t$ on $T^*M$. In the case where the considered K\"ahler
structures are Einstein they get several situations in which  the
parameters $a_1,\lambda $ are related by some algebraic relations.
In the general case, $(T^*M,G,J)$ has constant holomorphic
curvature.

 In this paper we study the singular case where the parameter
$a_1=\frac{A}{\lambda},A\in {\bf R}$. The class of the natural
almost complex structures $J$ on $T^*M$ that interchange the
vertical and horizontal distributions depends on two essential
parameters $\lambda$ and $b_1$. These parameters are smooth real
functions depending on the energy density $t$ on $T^*M$. From the
integrability condition for $J$ it follows that the base manifold
$M$ must have constant curvature $c$ and the second parameter
$b_1$ must be expressed as a rational function depending on the
first parameter $\lambda$ and its derivative. Of course, in the
obtained formula there are involved the constant $c$ and the
energy density $t$.

 A class of natural Riemannian metrics $G$ of diagonal type on $T^*M$
is defined by four parameters $c_1, c_2, d_1, d_2$ which are
smooth functions of~$t$. From the condition for $G$ to be
Hermitian with respect to $J$ we get two sets of proportionality
relations, from which we can get the parameters $c_1, c_2, d_1,
d_2$ as functions depending on one new parameter $\mu$ and the
parameter $\lambda $ involved in the expression of $J$.

 In the case where the fundamental $2$-form
$\phi$, associated to the class of complex structures $(G,J)$ is
closed, one finds that $\mu = \lambda^\prime$.

 Thus, we get
a class of K\"ahler structures $(G,J)$ on $T^*M$, depending on one
essential parameter $\lambda$.

 Finally, we prove that the obtained class of K\"ahler structures
on $T^*M$ is locally symmetric, Einstein and cannot have constant
holomorphic sectional curvature.

The manifolds, tensor fields and geometric objects we consider in
this paper, are assumed to be differentiable of class $C^{\infty}$
(i.e. smooth). We use the computations in local coordinates but
many results from this paper may be expressed in an invariant
form. The well known summation convention is used throughout this
paper, the range for the indices $h,i,j,k,l,r,s$ being
always${\{}1,...,n{\}}$ (see \cite{GheOpr}, \cite{OprPap}). We
shall denote by ${\Gamma}(T^*M)$ the module of smooth vector
fields on $T^*M$.

\vskip5mm {\large \bf 1. Some geometric properties of $T^*M$}
\vskip3mm

Let $(M,g)$ be a smooth $n$-dimensional Riemannian manifold and
denote its cotangent bundle by $\pi :T^*M\longrightarrow M$.
Recall that there is a structure of a $2n$-dimensional smooth
manifold on $T^*M$, induced from the structure of smooth
$n$-dimensional manifold  of $M$. From every local chart
$(U,\varphi )=(U,x^1,\dots ,x^n)$  on $M$, it is induced a local
chart $(\pi^{-1}(U),\Phi )=(\pi^{-1}(U),q^1,\dots , q^n,$
$p_1,\dots ,p_n)$, on $T^*M$, as follows. For a cotangent vector
$p\in \pi^{-1}(U)\subset T^*M$, the first $n$ local coordinates
$q^1,\dots ,q^n$ are  the local coordinates $x^1,\dots ,x^n$ of
its base point $x=\pi(p)$ in the local chart $(U,\varphi )$ (in
fact we have $q^i=\pi ^* x^i=x^i\circ \pi, \ i=1,\dots n)$. The
last $n$ local coordinates $p_1,\dots ,p_n$ of $p\in \pi^{-1}(U)$
are the vector space coordinates of $p$ with respect to the
natural basis $(dx^1_{\pi(p)},\dots , dx^n_{\pi(p)})$, defined by
the local chart $(U,\varphi )$,\ i.e. $p=p_idx^i_{\pi(p)}$.

 An $M$-tensor field of type $(r,s)$ on $T^*M$ is defined by sets
of $n^{r+s}$ components (functions depending on $q^i$ and $p_i$),
with $r$ upper indices and $s$ lower indices, assigned to induced
local charts $(\pi^{-1}(U),\Phi )$ on $T^*M$, such that the local
coordinate change rule is that of the local coordinate components
of a tensor field of type $(r,s)$ on the base manifold $M$ (see
\cite{Mok} for further details in the case of the tangent bundle).
An usual tensor field of type $(r,s)$ on $M$ may be thought of as
an $M$-tensor field of type $(r,s)$ on $T^*M$. If the considered
tensor field on $M$ is covariant only, the corresponding
$M$-tensor field on $T^*M$ may be identified with the induced
(pullback by $\pi $) tensor field on $T^*M$.

Some useful $M$-tensor fields on $T^*M$ may be obtained as
follows. Let $u,v:[0,\infty ) \longrightarrow {\bf R}$ be a smooth
functions and let $\|p\|^2=g^{-1}_{\pi(p)}(p,p)$ be the square of
the norm of the cotangent vector $p\in \pi^{-1}(U)$ ($g^{-1}$ is
the tensor field of type (2,0) having the components $(g^{kl}(x))$
which are the entries of the inverse of the matrix $(g_{ij}(x))$
defined by the components of $g$ in the local chart $(U,\varphi
)$). The components $u(\|p\|^2)g_{ij}(\pi(p))$, $p_i$,
$v(\|p\|^2)p_ip_j $ define $M$-tensor fields of types $(0,2)$,
$(0,1)$, $(0,2)$ on $T^*M$, respectively. Similarly, the
components $u(\|p\|^2)g^{kl}(\pi(p))$, $g^{0i}=p_hg^{hi}$,
$v(\|p\|^2)g^{0k}g^{0l}$ define $M$-tensor fields of type $(2,0)$,
$(1,0)$, $(2,0)$ on $T^*M$, respectively. Of course, all the
components considered above are in the induced local chart
$(\pi^{-1}(U),\Phi)$.

 The Levi Civita connection $\dot \nabla $ of $g$ defines a direct
sum decomposition

\begin{equation}
TT^*M=VT^*M\oplus HT^*M.
\end{equation}
of the tangent bundle to $T^*M$ into vertical distributions
$VT^*M= {\rm Ker}\ \pi _*$ and the horizontal distribution
$HT^*M$.

 If $(\pi^{-1}(U),\Phi)=(\pi^{-1}(U),q^1,\dots ,q^n,p_1,\dots ,p_n)$
is a local chart on $T^*M$, induced from the local chart
$(U,\varphi )= (U,x^1,\dots ,x^n)$, the local vector fields
$\frac{\partial}{\partial p_1}, \dots , \frac{\partial}{\partial
p_n}$ on $\pi^{-1}(U)$ define a local frame for $VT^*M$ over $\pi
^{-1}(U)$ and the local vector fields $\frac{\delta}{\delta
q^1},\dots ,\frac{\delta}{\delta q^n}$ define a local frame for
$HT^*M$ over $\pi^{-1}(U)$, where
$$
\frac{\delta}{\delta q^i}=\frac{\partial}{\partial
q^i}+\Gamma^0_{ih} \frac{\partial}{\partial p_h},\ \ \ \Gamma
^0_{ih}=p_k\Gamma ^k_{ih}
 $$
and $\Gamma ^k_{ih}(\pi(p))$ are the Christoffel symbols of $g$.

The set of vector fields $(\frac{\partial}{\partial p_1},\dots
,\frac{\partial}{\partial p_n}, \frac{\delta}{\delta q^1},\dots
,\frac{\delta}{\delta q^n})$ defines a local frame on $T^*M$,
adapted to the direct sum decomposition (1).

We consider
\begin{equation}
t=\frac{1}{2}\|p\|^2=\frac{1}{2}g^{-1}_{\pi(p)}(p,p)=\frac{1}{2}g^{ik}(x)p_ip_k,
\ \ \ p\in \pi^{-1}(U)
\end{equation}
the energy density defined by $g$ in the cotangent vector $p$. We
have $t\in [0,\infty)$ for all $p\in T^*M$.

From now on we shall work in a fixed local chart $(U,\varphi)$ on
$M$ and in the induced local chart $(\pi^{-1}(U),\Phi)$ on $T^*M$.

Now we shall present the following auxiliary result. \vskip5mm

{\bf Lemma 1}. \it If ~$n>1$ and $u,v$ are smooth functions on
$T^*M$ such that
$$
u g_{ij}+v p_ip_j=0,\ p\in \pi^{-1}(U)
$$
on the domain of any induced local chart on $T^*M$, then $u=0,\
v=0$.\rm \vskip 0.5cm

The proof is obtained easily by transvecting the given relation
with the components $g^{ij}$ of the tensor field $g^{-1}$ and
$g^{0j}$.\vskip5mm

\bf Remark. \rm From the relations of the type
$$
u g^{ij}+v g^{0i}g^{0j}=0,\ p\in \pi^{-1}(U),
$$
$$
u\delta ^i_j+vg^{0i} p_j=0,\ p\in \pi^{-1}(U),
$$
it is obtained, in a similar way, $u=v=0$.

\vskip5mm {\large \bf 2. A class of natural complex structures of
diagonal type on $T^*M$} \vskip3mm

 Consider the real valued smooth functions $\lambda,a_1,a_2,b_1,b_2$
on energy density  $t \in [0,\infty)$. We define a class of
natural almost complex structures $J$ of diagonal type on $T^*M$ ,
expressed in the adapted local frame by

\begin{equation}
~~~~~~~~~J{\frac{\delta}{\delta
q^i}}=J^{(1)}_{ij}(p){\frac{\partial}{\partial p_j}},
~~~~~~~~~~~~~J{\frac{\partial}{\partial
p_i}}=-J_{(2)}^{ij}(p){\frac{\delta}{\delta q^j}}.
\end{equation}

where,
\begin{equation}
J^{(1)}_{ij}(p)=a_1(t)g_{ij} + b_1(t)p_ip_j, ~~~
J_{(2)}^{ij}(p)=a_2(t)g^{ij} + b_2(t)g^{0i}g^{0j},~~~A\in {\bf
R^*}.
\end{equation}

 In this paper we study the singular case where

\begin{equation}
a_1(t)=\frac{A}{\lambda(t)}.
\end{equation}

 The components $J^{(1)}_{ij},J_{(2)}^{ij}$ define symmetric
 $M$-tensor fields of types $(0,2),(2,0)$ on $T^*M$, respectively.
\vskip5mm
{\bf Proposition 2}. \it The operator $J$ defines an
almost complex structure on $T^*M$ if and only if

\begin{equation}
a_1a_2=1,~~~ (a_1+2tb_1)(a_2+2tb_2)=1.
\end{equation}
\vskip2mm

 Proof. \rm The relations are obtained easily from the property $J^2=-I$ of
$J$ and Lemma 1.

From the relations (5), (6) we can obtain the explicit expression
of the parameter $a_2, b_2$
\begin{equation}
a_2=\frac{\lambda}{A},~~~~~~~~b_2=-\frac{\lambda^2 b_1}{A(A+2t
\lambda b_1)}.
\end{equation}

The obtained class of almost complex structures defined by the
tensor field $J$ on $T^*M$ is called \it class of natural almost
complex structures of diagonal type, \rm obtained from the
Riemannian metric $g$, by using the parameters $\lambda,\ b_1$. We
use the word diagonal for these almost complex structures, since
the $2n\times 2n$-matrix associated to $J$, with respect to the
adapted local frame $(\frac{\delta}{\delta q^1},\dots
,\frac{\delta}{\delta q^n},\frac{\partial}{\partial p_1},\dots
,\frac{\partial}{\partial p_n})$ has two $n\times n$-blocks on the
second diagonal
\begin{displaymath}
J= \left(
\begin{array}{cc}
0 & -J_{(2)}^{ij} \\
J^{(1)}_{ij} & 0
\end{array}
\right).
\end{displaymath}

\bf Remark. \rm From the conditions (6) it follows that the
coefficients $a_1,a_2,\\a_1+2tb_1,a_2+2tb_2$ cannot vanish and
have the same sign.

We assume that
\begin{equation}
\lambda (t) > 0 ~~~\forall ~t\geq0,~~~A>0.
\end{equation}

 It follows that~~~$a_1+2tb_1>0,~~~a_2+2tb_2>0 ~~~\forall ~t\geq0$, i.e.
\begin{equation}
 A+2t\lambda b_1>0,~~~\forall ~t\geq0.
\end{equation}

Now we shall study the integrability of the class of natural
almost complex structures defined by $J$ on $T^*M$. To do this we
need the following well known formulas for the brackets of the
vector fields $\frac{\partial}{\partial p_i}, \frac{\delta}{\delta
q^i},~ i=1,...,n$
\begin{equation}
[\frac{\partial}{\partial p_i},\frac{\partial}{\partial
p_j}]=0,~~~[\frac{\partial}{\partial p_i},\frac{\delta}{\delta
q^j}]=\Gamma^i_{jk}\frac{\partial}{\partial p_k},~~~
[\frac{\delta}{\delta q^i},\frac{\delta}{\delta q^j}]
=R^0_{kij}\frac{\partial}{\partial p_k},
\end{equation}
where $R^h_{kij}(\pi(p))$ are the local coordinate components of
the curvature tensor field of $\dot \nabla$ on $M$ and
$R^0_{kij}(p)=p_hR^h_{kij}$ . Of course, the components
 $R^h_{kij}$, $R^0_{kij}$ define M-tensor fields of types
 (1,3), (0,3) on $T^*M$, respectively.\\

Recall that the Nijenhuis tensor field $N$ defined by $J$ is given
by
$$
N(X,Y)=[JX,JY]-J[JX,Y]-J[X,JY]-[X,Y],\ \ \forall\ \ X,Y \in \Gamma
(T^*M).
$$
Then, we have $\frac{\delta}{\delta q^k}t =0,\
\frac{\partial}{\partial p_k}t = g^{0k}$. The expressions for the
components of $N$ can be obtained by a quite long, straightforward
computation, as follows\\

 {\bf Theorem 3. } {\it The Nijenhuis tensor field of the almost
complex structure $J$ on $T^*M$ is given by}
$$
\left\{
\begin{array}{l}
N(\frac{\delta}{\delta q^i},\frac{\delta}{\delta
q^j})=-\{\frac{A}{\lambda^3}[b_1\lambda
(\lambda+2t\lambda^{\prime})+A \lambda^{\prime}](\delta^h_ig_{jk}-
\delta^h_jg_{ik})+R^h_{kij}\}p_h\frac{\partial}{\partial p_k},
\\ \mbox{ } \\
N(\frac{\delta}{\delta q^i},\frac{\partial}{\partial
p_j})=-J_{(2)}^{kl}J_{(2)}^{jr}\{\frac{A}{\lambda^3}[b_1\lambda
(\lambda+2t\lambda^{\prime})+A \lambda^{\prime}](\delta^h_ig_{rl}-
\delta^h_rg_{il})+R^h_{lir}\}p_h\frac{\delta}{\delta q^k},
\\ \mbox{ } \\
N(\frac{\partial}{\partial p_i},\frac{\partial}{\partial
p_j})=-J_{(2)}^{ir}J_{(2)}^{jl}\{\frac{A}{\lambda^3}[b_1\lambda
(\lambda+2t\lambda^{\prime})+A \lambda^{\prime}](\delta^h_lg_{rk}-
\delta^h_rg_{lk})+R^h_{klr}\}p_h\frac{\partial}{\partial p_k}.
\end{array}
\right.
$$

\vskip2mm {\bf Theorem 4. } {\it The almost complex structure $J$
on $T^*M$ is integrable if and only if $(M,g)$ has  constant
sectional curvature $c$ and the function $b_1$ is given by
\begin{equation}
b_1=-\frac{c \lambda^3+A^2 \lambda^{\prime}}{A \lambda (\lambda+2t
\lambda^{\prime})}.
\end{equation}

The parameter $\lambda$ must fulfill the conditions
\begin{equation}
\lambda>0,~~~\frac{A^2-2ct \lambda^2}{\lambda+2t \lambda^{\prime}}
> 0~~~ \forall ~t \geq 0,~~~A>0.
\end{equation}}

 {\it Proof. }
  From the condition $N=0$, one obtains
$$
\{\frac{A}{\lambda^3}[b_1\lambda (\lambda+2t\lambda^{\prime})+A
\lambda^{\prime}](\delta^h_ig_{jk}-
\delta^h_jg_{ik})+R^h_{kij}\}p_h = 0
$$

Differentiating with respect to $p_l$, taking ~$p_h=0~~\forall ~h
\in {\{}1,...,n{\}}$, it follows that the curvature tensor field
of $\dot \nabla$ has the expression

$$
R^l_{kij}=-\frac{A}{\lambda(0)^3}[b_1(0)\lambda(0)^2+A
\lambda^{\prime}(0)]({\delta}^l_ig_{jk}-{\delta}^l_jg_{ik}).
$$
Using by the Schur theorem(in the case where $M$ is connected and
$dim M \geq 3$) it follows that $(M,g)$ has the constant sectional
curvature $c=-\frac{A}{\lambda(0)^3}[b_1(0)\lambda(0)^2+A
\lambda^{\prime}(0)]$ . Then we obtain the expression (11) of
$b_1$.

Conversely, if $(M,g)$ has constant curvature $c$ and $b_1$ is
given by (11), it follows in a straightforward way that $N = 0$.

By using (8),(9),(11) we obtain the conditions (12).\\

 The class of natural complex structures $J$ of diagonal type
on $T^*M$ depends on one essential parameter $\lambda$. The
components of $J$ are given by
\begin{equation}
J^{(1)}_{ij}=\frac{A}{\lambda}g_{ij}-\frac{c \lambda^3+A^2
\lambda^{\prime}}{A \lambda (\lambda+2t
\lambda^{\prime})}p_ip_j,~~~~~
J_{(2)}^{ij}=\frac{\lambda}{A}g^{ij}+\frac{c \lambda^3+A^2
\lambda^{\prime}}{A(A^2-2ct \lambda^2)}g^{0i}g^{0j}.
\end{equation}

\vskip5mm {\large \bf 3. A class of natural Hermitian structures
on $T^*M$}\vskip3mm

Consider the following symmetric $M-$tensor fields on $T^*M$,
defined by the components
\begin{equation}
G^{(1)}_{ij}=c_1 g_{ij}+d_1p_ip_j,\ \ \
G_{(2)}^{ij}=c_2g^{ij}+d_2g^{0i} g^{0j},
\end{equation}
where $c_1,c_2,d_1,d_2$ are smooth functions depending on the
energy density $t\in [0,\infty)$.

Obviously, $G^{(1)}$ is of type $(0,2)$ and $G_{(2)}$ is of type
$(2,0)$. We shall assume that the matrices defined by $G^{(1)}$
and $G_{(2)}$ are positive definite. This happens if and only if
\begin{equation}
c_1>0,~~c_2>0,~~c_1+2td_1>0,~~c_2+2td_2>0.
\end{equation}
Then the following class of Riemannian metrics may be considered
on $T^*M$
\begin{equation}
G=G^{(1)}_{ij}dq^idq^j+G_{(2)}^{ij}Dp_iDp_j,
\end{equation}
where $Dp_i=dp_i-\Gamma^0_{ij}dq^j$ is the absolute (covariant)
differential of $p_i$ with respect to the Levi Civita connection
$\dot\nabla$ of $g$. Equivalently, we have
$$
G(\frac{\delta}{\delta q^i},\frac{\delta}{\delta
q^j})=G^{(1)}_{ij},~~G(\frac{\partial}{\partial
p_i},\frac{\partial}{\partial p_j})=G_{(2)}^{ij},~~
G(\frac{\partial}{\partial p_i},\frac{\delta}{\delta q^j})=
G(\frac{\delta}{\delta q^j},\frac{\partial}{\partial p_i})=0.
$$
Remark that $HT^*M,~VT^*M$ are orthogonal to each other with
respect to $G$, but the Riemannian metrics induced from $G$ on
$HT^*M,~VT^*M$ are not the same, so the considered metric $G$ on
$T^*M$ is not a metric of Sasaki type. The $2n\times 2n$-matrix
associated to $G$, with respect to the adapted local frame
$(\frac{\delta}{\delta q^1},\dots ,\frac{\delta}{\delta
q^n},\frac{\partial}{\partial p_1},\dots ,\frac{\partial}{\partial
p_n})$ has two $n\times n$-blocks on the first diagonal
\begin{displaymath}
G= \left(
\begin{array}{cc}
G^{(1)}_{ij} & 0  \\
0 & G_{(2)}^{ij}
\end{array}
\right).
\end{displaymath}

The class of Riemannian metrics $G$ is called a \it class of
natural lifts of diagonal type \rm of $g$. Remark also that the
system of 1-forms $(Dp_1,...,Dp_n, dq^1,...,dq^n)$ defines a local
frame on $T^*T^*M$, dual to the local frame
$(\frac{\partial}{\partial p_1},\dots ,\frac{\partial}{\partial
p_n}, \frac{\delta}{\delta q^1},\dots ,\frac{\delta}{\delta
q^n})$, on $TT^*M$ over $\pi^{-1}(U)$ adapted to the direct sum
decomposition (1).

We shall consider another two $M$-tensor fields $H_{(1)}, \
H^{(2)}$ on $T^*M$, defined by the components
$$
H_{(1)}^{jk}=\frac{1}{c_1}g^{jk}-\frac{d_1}{c_1(c_1+2td_1)}g^{0j}g^{0k},
$$
$$
H^{(2)}_{jk}=\frac{1}{c_2}g_{jk}-\frac{d_2}{c_2(c_2+2td_2)}p_jp_k.
$$

The components $H_{(1)}^{jk}$ define an $M$-tensor field of type
$(2,0)$ and the components $H^{(2)}_{jk}$ define an $M$-tensor
field of type $(0,2)$. Moreover, the matrices associated to
$H_{(1)}, \ H^{(2)}$ are the inverses of the matrices associated
to  $G^{(1)}$ and  $G_{(2)}$, respectively. Hence we have
$$
G^{(1)}_{ij}H_{(1)}^{jk} = \delta_i^k,\ \ G_{(2)}^{ij}H^{(2)}_{jk}
= \delta^i_k.
$$

Now, we shall be interested in the conditions under which the
class of the metrics $G$ is Hermitian with respect to the class of
the complex structures $J$, considered in the previous section,
i.e.
$$
G(JX,JY)=G(X,Y),
$$
for all vector fields $X,Y$ on $T^*M$.

Considering the coefficients of $g_{ij}, g^{ij}$ in the conditions
\begin{equation}
\left\{
\begin{array}{l}
G(J\frac{\delta}{\delta q^i},J\frac{\delta}{\delta q^j})=
G(\frac{\delta}{\delta q^i},\frac{\delta}{\delta q^j}),
\\ \mbox{ } \\
G(J\frac{\partial}{\partial p_i},J\frac{\partial}{\partial p_j})=
G(\frac{\partial}{\partial p_i},\frac{\partial}{\partial p_j}),
\end{array}
\right.
\end{equation}
we can express the parameters $c_1,c_2$  with the help of the
parameters $a_1, a_2$ and a proportionality factor which must be
$\lambda = \lambda (t)$ ( see \cite{OprPor2}). Then
\begin{equation}
c_1 = \lambda a_1 = A,~~~~~~~~\ c_2=\lambda a_2 =
\frac{\lambda^2}{A},
\end{equation}
where the coefficients $a_1,a_2$ are given by (5) and (7).

Next, considering the coefficients of $p_ip_j,\ g^{0i}g^{0j}$ in
the relations (17), we can express the parameters
$c_1+2td_1,c_2+2td_2$ with help of the parameters $a_1+2tb_1,
a_2+2tb_2$ and a proportionality factor $\lambda+2t\mu $

\begin{equation}
\left\{
\begin{array}{l}
c_1+2td_1=(\lambda+2t\mu )(a_1+2tb_1),
\\ \mbox{ } \\
c_2+2td_2=(\lambda+2t\mu )(a_2+2tb_2).
\end{array}
\right.
\end{equation}
Remark that~ $\lambda (t)+2t \mu (t) > 0~~\forall~ t\geq 0$. It is
much more convenient to consider the proportionality factor in
such a form in the expression of the parameters
$c_1+2td_1,~c_2+2td_2$. Using by the relations (5), (7), (11),
(18) we can obtain easily from (19) the explicit expressions of
the coefficients $d_1,d_2$
\begin{equation}
\left\{
\begin{array}{l}
d_1=\frac{-(c\lambda^3+A^2\lambda^{\prime})+\mu(A^2-2ct\lambda^2)}{A(\lambda+2t\lambda^{\prime})},
\\ \mbox{ } \\
d_2=\frac{\lambda(c\lambda^3+A^2\lambda^{\prime})+\mu
A^2(\lambda+2\lambda^{\prime}t)}{A(A^2-2ct\lambda^2)}.
\end{array}
\right.
\end{equation}
Hence we may state:

{\bf Theorem 5.} \it Let $J$ be the class of natural, complex
structure of diagonal type on $T^*M$, given by (3) and (13). Let
$G$ be the class of the natural Riemannian metrics of diagonal
type on $T^*M$, given by (14), (18), (20).

Then we obtain a class of Hermitian structures $(G,J)$ on $T^*M$,
depending on two essential parameters $\lambda$ and $\mu$, which
must fulfill the conditions
\begin{equation}
\lambda > 0,~~~ \frac{A^2-2ct \lambda^2}{\lambda+2t
\lambda^{\prime}}
> 0,~~~\lambda +2t\mu >0,~~~\forall ~t\geq 0,~~~A>0.
\end{equation} \rm \vskip5mm

\vskip5mm {\large \bf 4. A class of K\"ahler structures on $T^*M$}
\vskip3mm
 Consider now the two-form $\phi $ defined by the class of
Hermitian structures $(G,J)$ on $T^*M$
$$
\phi (X,Y)=G(X,JY),
$$
for all vector fields $X,Y$ on $T^*M$.

Using by the expression of $\phi$ and computing the values $\phi
(\frac{\partial}{\partial p_i},\frac{\partial}{\partial p_j}),
\phi(\frac{\delta}{\delta q^i},\frac{\delta}{\delta q^j}),\\
\phi(\frac{\partial}{\partial p_i},\frac{\delta}{\delta q^j})$, we
obtain. \vskip5mm

 \bf Proposition 6. \it The expression of the
$2$-form $\phi $ in a local adapted frame
$(\frac{\partial}{\partial p_1},\dots ,\frac{\partial}{\partial
p_n}, \frac{\delta}{\delta q^1},\dots ,\frac{\delta}{\delta q^n})$
on $T^*M$, is given by
$$
\phi (\frac{\partial}{\partial p_i},\frac{\partial}{\partial
p_j})=0,\ \phi(\frac{\delta}{\delta q^i},\frac{\delta}{\delta
q^j})=0,\ \phi(\frac{\partial}{\partial p_i},\frac{\delta}{\delta
q^j})= \lambda \delta^i_j+\mu g^{0i}p_j,
$$
or, equivalently
\begin{equation}
\phi =(\lambda \delta^i_j+\mu g^{0i}p_j)Dp_i\wedge dq^j.
\end{equation}
\rm

\bf Theorem 7. \it The class of Hermitian structures $(G,J)$ on
$T^*M$ is K\"ahler if and only if
$$
\mu=\lambda ^\prime .
$$
Proof. \rm The expressions of $d\lambda,\ d\mu, \ dg^{0i}$ and
$dDp_i$ are obtained in a straightforward way, by using the
property $\dot \nabla_k g_{ij}=0$ (hence $\dot \nabla_k g^{ij}=0$)
$$
d\lambda = \lambda ^\prime g^{0i}Dp_i,\ d\mu =\mu ^\prime
g^{0i}Dp_i,\ dg^{0i}=g^{ik}Dp_k-g^{0h}\Gamma ^i_{hk}dq^k,
$$
$$
dDp_i=-\frac{1}{2}R^0_{ikl}dq^k\wedge dq^l+ \Gamma
^l_{ik}dq^k\wedge Dp_l.
$$
Then we have
$$
d\phi =(d\lambda \delta^i_j+d\mu g^{0i}p_j+ \mu dg^{0i}p_j+\mu
g^{0i}dp_j)\wedge Dp_i\wedge dq^j+
$$
$$
+(\lambda \delta^i_j+\mu g^{0i}p_j)dDp_i\wedge dq^j.
$$
By replacing the expressions of $d\lambda , d\mu , dg^{0i}$ and
$d\dot\nabla y^h$, then using, again, the property $\dot\nabla
_kg_{ij}=0$, doing some algebraic computations with the exterior
products, then using the well known symmetry properties of
$g_{ij}, \Gamma ^h_{ij},$ and of the Riemann-Christoffel tensor
field, as well as the Bianchi identities, it follows that
$$
d\phi =\frac{1}{2}(\lambda ^\prime -\mu)g^{0h}Dp_h\wedge
Dp_i\wedge dq^i.
$$
Therefore we have $d\phi =0$ if and only if $\mu =\lambda ^\prime
$. \vskip5mm

 {\bf Remark.} The class of natural K\"ahler structures
of diagonal type defined by $(G,J)$ on $T^*M$ depends on one
essential parameter $\lambda$.

The paramater $\lambda$ must fulfill the conditions
\begin{equation}
\lambda > 0,~~~A^2-2ct \lambda^2>0,~~~\lambda +
2t\lambda^{\prime}>0,~~~\forall ~t \geq 0,~~~A>0.
\end{equation}

The components of the class of K\"ahler metrics $G$ are given by
\begin{equation}
\left\{
\begin{array}{l}
G^{(1)}_{ij} = Ag_{ij} - \frac{c\lambda^2}{A}p_ip_j,
\\ \mbox{ } \\
G^{ij}_{(2)} = \frac{\lambda^2}{A}g^{ij} + \frac{c\lambda^4 +
2A^2\lambda^{\prime}(\lambda + \lambda^{\prime}t)}{A(A^2 -
2ct\lambda^2)}g^{0i}g^{0j}.
\end{array}
\right.
\end{equation}

We obtain, too
\begin{equation}
\left\{
\begin{array}{l}
H^{jk}_{(1)}=\frac{1}{A}g^{jk}+\frac{c\lambda^2}{A(A^2-2ct\lambda^2)}g^{0j}g^{0k},
\\ \mbox{ } \\
H^{(2)}_{jk}=\frac{A}{\lambda^2}g_{jk}-\frac{c\lambda^4 +
2A^2\lambda^{\prime}(\lambda +
\lambda^{\prime}t)}{A\lambda^2(\lambda+2\lambda^{\prime}t)^2}p_jp_k.
\end{array}
\right.
\end{equation}
\vskip5mm

{\large \bf 5. A class of locally symmetric K\"ahler Einstein
structures on $T^*M$} \vskip3mm

The Levi Civita connection $\nabla$ of the Riemannian manifold
$(T^*M,G)$ is determined by the conditions
$$
\nabla G=0,~~~~~  T =0,
$$
where $T$ is its torsion tensor field. The explicit expression of
this connection is obtained from the formula
$$
2G({\nabla}_XY,Z)=X(G(Y,Z))+Y(G(X,Z))-Z(G(X,Y))+
$$
$$
+G([X,Y],Z)-G([X,Z],Y)-G([Y,Z],X); ~~~~~~ \forall\
X,Y,Z~{\in}~{\Gamma}(T^*M).
$$

The final result can be stated as follows. \\

\bf Theorem 8. {\it The Levi Civita connection ${\nabla}$ of $G$
has the following expression in the local adapted frame
$(\frac{\delta}{\delta q^1},\dots ,\frac{\delta}{\delta
q^n},\frac{\partial}{\partial p_1},\dots ,\frac{\partial}{\partial
p_n}):$
\begin{equation}
\left\{
\begin{array}{l}
 \nabla_\frac{\partial}{\partial
p_i}\frac{\partial}{\partial p_j}
=Q^{ij}_h\frac{\partial}{\partial p_h},\ \ \ \ \ \
\nabla_\frac{\delta}{\delta q^i}\frac{\partial}{\partial
p_j}=-\Gamma^j_{ih}\frac{\partial}{\partial
p_h}+P^{hj}_i\frac{\delta}{\delta q^h},
\\ \mbox{ } \\
\nabla_\frac{\partial}{\partial p_i}\frac{\delta}{\delta
q^j}=P^{hi}_j\frac{\delta}{\delta q^h},\ \ \ \ \ \
\nabla_\frac{\delta}{\delta q^i}\frac{\delta}{\delta
q^j}=\Gamma^h_{ij}\frac{\delta}{\delta
q^h}+S_{hij}\frac{\partial}{\partial p_h},
\end{array}
\right.
\end{equation}

where $Q^{ij}_h, P^{hi}_j, S_{hij}$ are $M$-tensor fields on
$T^*M$, defined by
\begin{equation}
\left\{
\begin{array}{l}
Q^{ij}_h = \frac{1}{2}H^{(2)}_{hk}(\frac{\partial}{\partial
p_i}G_{(2)}^{jk}+ \frac{\partial}{\partial p_j}G_{(2)}^{ik}
-\frac{\partial}{\partial p_k}G_{(2)}^{ij}),
\\ \mbox{ } \\
P^{hi}_j=\frac{1}{2}H_{(1)}^{hk}(\frac{\partial}{\partial
p_i}G^{(1)}_{jk}-G_{(2)}^{il}R^0_{ljk}),
\\ \mbox{ } \\
S_{hij}=-\frac{1}{2}H^{(2)}_{hk}\frac{\partial}{\partial
p_k}G^{(1)}_{ij}+\frac{1}{2}R^0_{hij}.
\end{array}
\right.
\end{equation}

 \rm

Assuming that the base manifold $(M,g)$ has constant sectional
curvature $c$ and replacing the expressions of the involved
$M$-tensor fields, one obtains

\begin{equation}
\left\{
\begin{array}{l}
Q^{ij}_h = \frac{c \lambda^3+A^2 \lambda^{\prime}}{A \lambda
(\lambda+2t \lambda^{\prime})}J^{ij}_{(2)}p_h +
\frac{\lambda^{\prime}}{\lambda^2}(\delta^i_hg^{0j}+\delta^j_hg^{0i})
+ \frac{\lambda \lambda^{\prime \prime } - 3
\lambda^{\prime^2}}{\lambda(\lambda+2\lambda^{\prime}t)}g^{0i}g^{0j}p_h,
\\ \mbox{ } \\
P^{hi}_j=-\frac{c\lambda}{A}J^{hi}_{(2)}p_j,
\\ \mbox{ } \\
S_{hij}=\frac{c\lambda}{A}J^{(1)}_{hj}p_i.
\end{array}
\right.
\end{equation}

The curvature tensor field $K$ of the connection $\nabla $ is
obtained from the well known formula
$$
K(X,Y)Z=\nabla_X\nabla_YZ-\nabla_Y\nabla_XZ-\nabla_{[X,Y]}Z,\ \ \
\ \forall\ X,Y,Z\in \Gamma (T^*M).
$$

The components of curvature tensor field $K$ with respect to the
adapted local frame $(\frac{\delta}{\delta q^1},\dots
,\frac{\delta}{\delta q^n},\frac{\partial}{\partial p_1},\dots
,\frac{\partial}{\partial p_n})$ are obtained easily:
\begin{equation}
\left\{
\begin{array}{l}
K(\frac{\delta}{\delta q^i},\frac{\delta}{\delta
q^j})\frac{\delta}{\delta
q^k}=\frac{c}{A}(\delta^h_iG^{(1)}_{jk}-\delta
^h_jG^{(1)}_{ik})\frac{\delta}{\delta q^h},
\\ \mbox{ } \\
K(\frac{\delta}{\delta q^i},\frac{\delta}{\delta
q^j})\frac{\partial}{\partial
p_k}=\frac{c}{A}G^{sk}_{(2)}(J^{(1)}_{js}J^{(1)}_{ih}-J^{(1)}_{is}J^{(1)}_{jh})\frac{\partial}{\partial
p_h},
\\ \mbox{ } \\
K(\frac{\partial}{\partial p_i},\frac{\partial}{\partial
p_j})\frac{\delta}{\delta
q^k}=\frac{c}{A}G^{(1)}_{sk}(J^{js}_{(2)}J^{ih}_{(2)}-J^{is}_{(2)}J^{jh}_{(2)})\frac{\delta}{\delta
q^h},
\\ \mbox{ } \\
K(\frac{\partial}{\partial p_i},\frac{\partial}{\partial
p_j})\frac{\partial}{\partial p_k} =\frac{c}{A}(\delta
^i_hG^{jk}_{(2)}-\delta ^j_hG^{ik}_{(2)})\frac{\partial}{\partial
p_h},
\\ \mbox{ } \\
K(\frac{\partial}{\partial p_i},\frac{\delta}{\delta
q^j})\frac{\delta}{\delta
q^k}=\frac{c}{A}J^{(1)}_{kh}\frac{\partial}{\partial p_i}(\lambda
p_j)\frac{\partial}{\partial p_h},
\\ \mbox{ } \\
K(\frac{\partial}{\partial p_i},\frac{\delta}{\delta
q^j})\frac{\partial}{\partial
p_k}=-\frac{c}{A}J^{(kh)}_{(2)}\frac{\partial}{\partial
p_i}(\lambda p_j)\frac{\delta}{\delta q^h}.
\end{array}
\right.
\end{equation}
{\bf Remark.} From the local coordinates expression of the
curvature tensor field K we obtain that the class of K\"ahler
structures $(G,J)$ on $T^*M$ cannot have constant holomorphic
sectional curvature.\\

The Ricci tensor field Ric of $\nabla$ is defined by the formula:
$$
Ric(Y,Z)=trace(X\longrightarrow K(X,Y)Z),\ \ \ \forall\  X,Y,Z\in
\Gamma (T^*M).
$$
It follows
$$
\left\{
\begin{array}{l}
Ric(\frac{\delta}{\delta q^i},\frac{\delta}{\delta
q^j})=\frac{cn}{A}G^{(1)}_{ij},
\\ \mbox{ } \\
 Ric(\frac{\partial}{\partial
p_i},\frac{\partial}{\partial p_j})=\frac{cn}{A}G^{ij}_{(2)},
\\ \mbox{ } \\
Ric(\frac{\partial}{\partial p_i},\frac{\delta}{\delta
q^j})=Ric(\frac{\delta}{\delta q^j},\frac{\partial}{\partial
p_i})=0.
\end{array}
\right.
$$
Thus
\begin{equation}
 Ric=\frac{cn}{A}G.
\end{equation}

By using the expressions (29) and the relations $\nabla G =
0,~\nabla J = 0 $ we have computed the covariant derivatives of
curvature tensor field K in the local adapted frame
$(\frac{\delta}{\delta q^i},\frac{\partial}{\partial p_i})$ with
respect to the connection $\nabla$ and we obtained in all twelve
cases the result is zero.

 Hence we may state our main
result.\\

{\bf Theorem 9.} {\it Assume that the Riemannian manifold $(M,g)$
has constant sectional curvature $c$. Let $J$ be the class of
natural, complex structure of diagonal type on $T^*M$, given by
(3) and (13). Let $G$ be the class of the natural Riemannian
metrics of diagonal type on $T^*M$, given by (14) and (24).

Then $(G,J)$ is a class of locally symmetric K\"ahler Einstein
structures on $T^*M$, depending on one essential parameter
$\lambda$, which must fulfill the conditions (23):
$$
 \lambda
> 0,~~~A^2-2ct \lambda^2>0,~~~\lambda +
2t\lambda^{\prime}>0,~~~\forall~ t \geq 0,~~~A>0.
$$
\rm

{\bf Remark.} In the case where $c<0$, the function $\lambda = t^m
+B,~m,B\in{\bf R_+}$, fulfill the conditions (23).

In the case where $c>0$, the function $\lambda =
\frac{A}{\sqrt{2ct} + B},~B \in {\bf R_+}$, fulfill the conditions
(23).

In the case where $c=0$, we obtain that $(M,g)$ is a flat manifold
and $G$ is a class of flat metrics on $T^*M$.

\end{document}